\documentclass[12pt]{amsart}
\usepackage{amssymb}

\usepackage[latin1]{inputenc}
\usepackage{mathptmx}

\newcommand{\Z}{\mathbb{Z}}

\newcommand{\N}{\mathbb{N}}

\newcommand{\ab}{{\bf a}}

\DeclareMathOperator{\pnt}{\raise 0.5mm \hbox{\large\bf.}}

\DeclareMathOperator{\codim}{codim}

\DeclareMathOperator{\depth}{depth}

\DeclareMathOperator{\gin}{gin}
\DeclareMathOperator{\chara}{char}
\DeclareMathOperator{\Tor}{Tor}
\DeclareMathOperator{\projdim}{proj\,dim}
\DeclareMathOperator{\reg}{reg}

\def\+#1{\relax\ifmmode\if\noexpand #1\relax \mathop{\kern
    0pt^+{#1}}\nolimits\else \kern 0pt^+\!#1 \fi\else$^*$#1\fi}

\newtheorem{thm}{\bf Theorem}[section]
\newtheorem{lem}[thm]{\bf Lemma}
\newtheorem{cor}[thm]{\bf Corollary}

\theoremstyle{definition}

\newtheorem{rem}[thm]{\bf Remark}
\newtheorem{ex}[thm]{\bf Example}

\theoremstyle{plain}
\newtheorem*{thm*}{Theorem}

\title{Betti numbers and shifts in minimal graded free resolutions}
\subjclass{05E99; 13C14; 13D02}
\author{Tim R\"omer}
\address{Universit\"at Osnabr\"uck, Institut f\"ur Mathematik, 49069 Osnabr\"uck, Germany}
\email{troemer@uos.de}

\begin{document}

\begin{abstract}
Let $S=K[x_1,\ldots,x_n]$ be a polynomial ring  and $R=S/I$
where $I \subset S$ is a graded ideal.
The Multiplicity Conjecture of Herzog, Huneke, and Srinivasan
which was recently proved using the Boij--S\"oderberg theory states that
the multiplicity of $R$ is bounded above by a function of the maximal shifts in the minimal graded free resolution of $R$ over $S$ as well as bounded below by a function of the minimal shifts if $R$ is Cohen--Macaulay. In this paper we study the related problem to show that the total Betti-numbers of $R$ are also bounded above by a function of the shifts in the minimal graded free resolution of $R$ as well as bounded below by another function of the shifts if $R$ is Cohen--Macaulay. We also discuss the cases when these bounds are sharp.
\end{abstract}

\maketitle

\section{Introduction}
\label{intro}
Let $S=K[x_1,\ldots,x_n]$ be a polynomial ring over a field $K$
equipped with the standard grading  by setting $\deg(x_i)=1$.
We consider a standard graded $K$-algebra $R=S/I$ where $I\subset S$ is a graded ideal
and the minimal graded free resolution of $R$:
$$
0
\to
\bigoplus_{j\in \Z} S(-j)^{\beta^S_{p,j}(R)}
\to
\cdots
\to
\bigoplus_{j\in \Z} S(-j)^{\beta^S_{1,j}(R)}
\to S \to 0
$$
where $\beta^S_{i,j}(R)=\dim_K \Tor^S_i(R,K)_{j}$  are the graded Betti numbers and $p=\projdim(R)$
is the projective dimension of $R$. Let $\beta^S_{i}(R)=\sum_{j\in\Z}\beta^S_{i,j}(R)$ be the $i$-th total Betti number of $R$.
Recall that $R$ has a \emph{pure resolution} if the resolution has the following shape:
$$
0
\to
S(-d_p)^{\beta^S_{p}(R)}
\to
\cdots
\to
S(-d_1)^{\beta^S_{1}(R)}
\to S \to 0
$$
for some numbers $d_1,\ldots,d_p$.
Let $e(R)$ denote the multiplicity of $R$.
If $R$ is Cohen--Macaulay with a pure resolution,
then
Herzog and K\"uhl \cite{HEKU84}, and Huneke and Miller  \cite{HUMI85}
observed that
the following formulas hold:
$$
e(R) = \frac{1}{p!}\prod_{i=1}^{p} d_i
\text{ and }
\beta_i^S(R)=(-1)^{i+1} \prod_{j\neq i}\frac{d_j}{d_j-d_i}
\text{ for } i=1,\ldots,p.
$$
Consider for $1\leq i \leq p$ the numbers
$$
M_i=\max\{j \in \Z \colon \beta^S_{i,j}(R)\neq 0\}
\text{ and }
m_i=\min\{j \in \Z \colon \beta^S_{i,j}(R)\neq 0\}.
$$
In the last years many people studied the Multiplicity Conjecture of
Herzog, Huneke and Srinivasan (see \cite{HESR98} and \cite{HUMI85})
which states in its original form that if $R=S/I$ is  Cohen--Macaulay, then
$$
\frac{1}{p!}\prod_{i=1}^{p} m_i \leq e(R) \leq \frac{1}{p!}\prod_{i=1}^{p} M_i.
$$
Migliore, Nagel and the author \cite{MINARO05} extended this conjecture by the questions that
we have equality below or above if and only if $R$ has a pure resolution.
This conjecture is proved as a corollary of the Boij--S\"oderberg theory which was conjectured and developed partly by Boij--S\"oderberg \cite{BOSO} and then completely  proved by Eisenbud--Schreyer \cite{EISCH} (see also \cite{BOSO} and \cite{EIFLWE}). A natural question is whether under the Cohen--Macaulay assumption
the $i$-th total Betti number $\beta^S_{i}(R)$
can also be bounded by using the shifts in the minimal graded free resolution of $R$.
A natural guess for bounds is
\begin{eqnarray}
\label{guess1}
\prod_{1\leq j< i}
\frac{m_j}{m_i-m_j}
\cdot
\prod_{i<j\leq p}
\frac{m_j}{m_j-m_i}
\leq
\beta_i^S(R)
\leq
\prod_{1\leq j<i}
\frac{M_j}{M_i-M_j}
\cdot
\prod_{i<j\leq p}
\frac{M_j}{M_j-M_i}
\end{eqnarray}
for $i=1,\dots,p$.
We show that these bounds hold if $R$ is a complete intersection
and if $I$ is componentwise linear.
Moreover, in these cases we have equality above or below for all $i$ if and only if $R$ has a pure resolution.
In general these bounds are not valid. Indeed, we give a counterexample in Example \ref{counterguess1}.
For Cohen--Macaulay algebras with \emph{strictly quasi-pure resolutions},
i.e.\ $m_i>M_{i-1}$ for all $i$,
we show the bounds
\begin{eqnarray}
\label{guess2}
\prod_{1\leq j< i}
\frac{m_j}{M_i-m_j}
\cdot
\prod_{i<j\leq p}
\frac{m_j}{M_j-m_i}
\leq
\beta_i^S(R)
\leq
\prod_{1\leq j< i}
\frac{M_j}{m_i-M_j}
\cdot
\prod_{i<j\leq p}
\frac{M_j}{m_j-M_i}
\end{eqnarray}
for $i=1,\dots,p$.
Again we have equality below or above for all $i$ if and only if $R$ has a pure resolution.
Observe that
$$
\prod_{1\leq j< i}
\frac{m_j}{M_i-m_j}
\cdot
\prod_{i<j\leq p}
\frac{m_j}{M_j-m_i}
\leq
\prod_{1\leq j< i}
\frac{m_j}{m_i-m_j}
\cdot
\prod_{i<j\leq p}
\frac{m_j}{m_j-m_i},
$$
because
$M_i-m_j \geq m_i-m_j>0$ for $1\leq j< i$
and
$M_j-m_i\geq m_j-m_i> 0$ for $i<j\leq p$ respectively.
Thus the weaker lower bounds in (\ref{guess2}) hold also for all cases where the lower bounds in (\ref{guess1}) are valid.
But the numbers $\prod_{1\leq j< i}
\frac{M_j}{m_i-M_j}
\cdot
\prod_{i<j\leq p}
\frac{M_j}{m_j-M_i}
$ may be negative  and thus are not candidates for upper bounds in general.
Note that the Cohen--Macaulay assumption for the lower bound (\ref{guess2}) is essential.
We construct a non Cohen--Macaulay ideal as a counterexample in Example \ref{counterguess2}.
We have that
$$
\prod_{1\leq j<i}
\frac{M_j}{M_i-M_j}
\cdot
\prod_{i<j\leq p}
\frac{M_j}{M_j-M_i}
\leq
\frac{1}{(i-1)!\cdot (p-i)!}
\prod_{j\neq i}  M_j
$$
because in the Cohen--Macaulay case we have
$M_i-M_j \geq i-j$ for $1\leq j<i$
and
$M_j-M_i \geq j-i$ for $i<j\leq p$
respectively.
Hence one might still ask if the upper bound
\begin{eqnarray}
\label{guess3}
\beta_i^S(R)
\leq
\frac{1}{(i-1)!\cdot (p-i)!}
\prod_{j\neq i}  M_j
\end{eqnarray}
is valid for $i=1,\dots,p$. In addition to the cases that the bounds in (\ref{guess1})
hold if $R$ is a complete intersection and if $I$ is componentwise linear,
the bounds in (\ref{guess2}) hold if $R$ has a strictly quasi-pure resolution,
using the Boij--S\"oderberg theory we show that the lower bounds in (\ref{guess2}) and the upper bounds in (\ref{guess3}) hold if $S/I$ is Cohen--Macaulay. Moreover, we discuss the case where we have equality everywhere. See also \cite{MIR07} for related results.  Some remarks on possible upper bounds for non Cohen--Macaulay algebras are included in this paper.

We are grateful to Prof.\ J. Herzog
for inspiring discussions on the subject of this paper.
\section{Complete Intersections}
\label{completeintersections}
One of the first examples of Cohen--Macaulay algebras are complete
intersection.
For this we consider a complete intersection
$R=S/I$ where $I=(f_1,\dots,f_p)$ is a graded ideal
generated by a regular sequence
$f_1,\dots,f_p$.
Let $\deg(f_i)=d_i$ for $i=1,\dots,p$.
Without loss of generality we assume that
$d_1\geq \dots \geq d_p$.
The Koszul complex gives rise to a minimal graded free resolution of $R$
and thus we get that
\begin{eqnarray*}
\beta_i(R) &=& \binom{p}{i},\\
M_i        &=& d_1+\dots+d_i,\\
m_i        &=& d_p+\dots+d_{p-i+1}
\end{eqnarray*}
for $i=1,\dots,p$.
Note that $R$ has a pure resolution if and only if $d_1=\dots=d_p$.
The ideal $I$ has a linear resolution if and only if $d_1=\dots=d_p=1$.
Using these facts we prove:
\begin{thm}
\label{thmci}
Let $R=S/I$ be a complete intersection as described above.
Then:
\begin{enumerate}
 \item
We have for $i=1,\dots,p$ that
$$
\beta_i^S(R)
\leq
\prod_{1\leq j<i}
\frac{M_j}{M_i-M_j}
\cdot
\prod_{i<j\leq p}
\frac{M_j}{M_j-M_i}
\leq
\frac{1}{(i-1)!\cdot (p-i)!}
\prod_{j\neq i}  M_j.
$$
The first upper bound is reached for all $i$ if and only if $R$ has a pure re\-so\-lu\-tion.
Every upper bound is reached for all $i$ if and only if $I$ has a linear resolution.
\item
We have for $i=1,\dots,p$ that
$$
\beta_i^S(R)
\geq
\prod_{1\leq j< i}
\frac{m_j}{m_i-m_j}
\cdot
\prod_{i<j\leq p}
\frac{m_j}{m_j-m_i}
\geq
\prod_{1\leq j< i}
\frac{m_j}{M_i-m_j}
\cdot
\prod_{i<j\leq p}
\frac{m_j}{M_j-m_i}.
$$
Every lower bound is reached for all $i$ if and only if
$R$ has a pure resolution.
\end{enumerate}
\end{thm}
\begin{proof}
(i):
To prove the upper bound we compute
for $p\geq j>i$ that
$$
\frac{M_j}{M_j-M_i}
=
\frac{d_1+\dots+d_j}{d_{i+1}+\dots+d_j}
=
\frac{d_1+\dots+d_i}{d_{i+1}+\dots+d_j} +1
\geq
\frac{i\cdot d_i}{(j-i)\cdot d_{i+1}} +1
\geq
\frac{i}{j-i} +1
=
\frac{j}{j-i}
$$
and for $1\leq j<i$ that
$$
\frac{M_j}{M_i-M_j}
=
\frac{d_1+\dots+d_j}{d_{j+1}+\dots+d_i}
\geq
\frac{j \cdot d_j}{(i-j)\cdot d_{j+1}}
\geq
\frac{j}{(i-j)}.
$$
Observe that we have equality for all integers $i,j$
if and only if $d_1=\dots=d_p$.
Thus
\begin{eqnarray*}
\beta_i^S(R)
&=&
\binom{p}{i}  \\
&=&
\frac{i-1}{1}
\cdot
\frac{i-2}{2}
\cdots
\frac{1}{i-1}
\cdot
\frac{p}{p-i}
\cdot
\frac{p-1}{p-1-i}
\cdots
\frac{i+1}{1}
\\
& \leq &
\frac{M_{i-1}}{M_i-M_{i-1}}
\cdot
\frac{M_{i-2}}{M_i-M_{i-2}}
\cdots
\frac{M_1}{M_i-M_1}
\cdot
\frac{M_p}{M_p-M_i}
\cdot
\frac{M_{p-1}}{M_{p-1}-M_i}
\cdots
\frac{M_{i+1}}{M_{i+1}-M_i}
\\
&= &
\prod_{1\leq j<i}
\frac{M_j}{M_i-M_j}
\cdot
\prod_{i<j\leq p}
\frac{M_j}{M_j-M_i}
\\
&\leq&
\frac{1}{(i-1)!\cdot (p-i)!}
\prod_{j\neq i}  M_j
\end{eqnarray*}
where the last inequality was observed in Section \ref{intro}.
Moreover, we have that
$\beta_i^S(R)
=\prod_{1\leq j<i}
\frac{M_j}{M_i-M_j}
\cdot
\prod_{i<j\leq p}
\frac{M_j}{M_j-M_i}$
for all $1\leq i \leq p$
if and only if $R$ has a pure resolution.
It is also easy to see that $\beta_i^S(R)=\frac{1}{(i-1)!\cdot (p-i)!}
\prod_{j\neq i}  M_j$
for all $1\leq i \leq p$
if and only if $I$ has a linear resolution.

(ii):
Similarly, it follows from
$$
\frac{m_j}{m_j-m_i}
=
\frac{d_p+\dots+d_{p-j+1}}{d_{p-i}+\dots+d_{p-j+1}}
=
\frac{d_p+\dots+d_{p-i+1}}{d_{p-i}+\dots+d_{p-j+1}} +1
$$
$$\leq
\frac{i \cdot d_{p-i+1} }{(j-i)\cdot d_{p-i}} +1
\leq
\frac{i}{(j-i)} +1
=
\frac{j}{(j-i)}
$$
for $p\geq j>i$
and
$$
\frac{m_j}{m_i-m_j}
=
\frac{d_p+\dots+d_{p-j+1}}{d_{p-j}+\dots+d_{p-i+1}}
\leq
\frac{j \cdot d_{p-j+1}}{(i-j)\cdot d_{p-j}}
\leq
\frac{j}{(i-j)}
$$
for $1\leq j<i$ that
\begin{eqnarray*}
\beta_i^S(R)
&=&
\binom{p}{i}  \\
&\geq&
\prod_{1\leq j< i}
\frac{m_j}{m_i-m_j}
\cdot
\prod_{i<j\leq p}
\frac{m_j}{m_j-m_i}
\\
&\geq&
\prod_{1\leq j< i}
\frac{m_j}{M_i-m_j}
\cdot
\prod_{i<j\leq p}
\frac{m_j}{M_j-m_i}.
\end{eqnarray*}
The last inequality was observed in Section \ref{intro}.
Again we have equations everywhere for all $1\leq i \leq p$
if and only if $R$ has a pure resolution.

This concludes the proof.
\end{proof}

\begin{rem}
Instead of this direct approach one can also use the Boij--S\"oderberg theory (see \cite{BOSO}, \cite{BOSO2}, \cite{EIFLWE}and \cite{EISCH}). See Section \ref{boijsoederberg} for details
where we obtain beside other things again the lower bounds in (\ref{guess2}) and the upper bounds in (\ref{guess3}) using this approach.
\end{rem}

\section{Ideals with strictly quasi-pure resolutions}
\label{section_quasipure}

Motivated by the results of Section \ref{completeintersections}
one could hope that the bounds in (\ref{guess1}) are always valid.
This is not the case as the following example shows.

\begin{ex}
\label{counterguess1}
We consider the following situation.
Let $S=K[x_1,\dots,x_6]$ be a polynomial ring in
$6$ variables and consider the graded ideal
$I=(x_1x_2, x_1x_3, x_2x_4-x_5x_6, x_3x_4)$.
Using for example CoCoA \cite{COCOA}
one checks that $S/I$ is Cohen--Macaulay of dimension 3
and it has the minimal graded free resolution:
$$
0 \to S^2(-5) \to S^2(-3)\oplus S^3(-4) \to S^4(-2) \to S \to 0
$$
which is not pure.
We have
$$
M_1=m_1=2,\
M_2=4, m_2=3,\
M_3=m_3=5.
$$
But
$$
\frac{M_2}{M_2-M_1}\cdot \frac{M_3}{M_3-M_1}
=
\frac{4}{2}\cdot \frac{5}{3}
=
\frac{20}{6}
<
4=\beta^S_1(R)
$$
and hence the upper bound of (\ref{guess1}) is not valid.
Moreover,
$$
\frac{m_2}{m_2-m_1}\cdot \frac{m_3}{m_3-m_1}
=
\frac{3}{1}\cdot \frac{5}{3}
=
5 >  4=\beta^S_1(R).
$$
Thus also the lower bound of (\ref{guess1}) is false in general.
But the resolution is strictly quasi-pure
since $m_i>M_{i-1}$ for all $1\leq i \leq 3$.
Note that the bounds in (\ref{guess2}) hold.
Indeed, e.g.\ for $\beta^S_1(R)$ we have
$$
\frac{M_2}{m_2-M_1}\cdot \frac{M_3}{m_3-M_1}
=
\frac{4}{1}\cdot \frac{5}{3}
=
\frac{20}{3}
>4=\beta^S_1(R)
$$
and
$$
\frac{m_2}{M_2-m_1}\cdot \frac{m_3}{M_3-m_1}
=
\frac{3}{2}\cdot \frac{5}{3}
=
\frac{15}{6}
<  4=\beta^S_1(R).
$$
\end{ex}

We recall the following well-known result
which is due to Peskine and Szpiro \cite{PESZ}.

\begin{lem}
\label{quasihelpers}
Let $I \subset S$ be a graded ideal such that $R=S/I$
is Cohen--Macaulay and let $p=\projdim(R)$.
Then:
\begin{enumerate}
\item
$\sum_{i=1}^p (-1)^i \sum_{j} \beta^S_{ij}(R)=\sum_{i=1}^p (-1)^i \beta^S_{i}(R) = -1 $.
\item
$
\sum_{i=1}^p (-1)^i \sum_{j} j^k \cdot \beta^S_{ij}(R)
=
0$ for  $1\leq k \leq p-1$.
\end{enumerate}

\end{lem}
\begin{proof}
We have
$\sum_{i=1}^p (-1)^i \sum_{j} \beta^S_{ij}(R)=\sum_{i=1}^p (-1)^i \beta^S_{i}(R) = -\beta_0^S(R)=-1$.
For a proof of the other equalities see also e.g.\ \cite[Lemma 1.1]{HESR98}.
\end{proof}
We see that the graded Betti numbers satisfy a certain system of equations which sometimes is nowadays also called the \emph{Herzog-K\"uhl equations}.
Note that if $R$ has a pure resolution, then using this system, Cramer's rule and the Vandermonde
determinant it is not difficult to prove
the formulas of the multiplicity and the total Betti-numbers in \cite{HEKU84} and  \cite{HUMI85}.
Recall from \cite{HESR98} that $R$ has a \emph{quasi-pure resolution}
if $m_i\geq M_{i-1}$ for all $i$. Unfortunately, we can not prove in general
the bounds in (\ref{guess2}) for the total Betti-numbers in this case.
We say that $R$ has a \emph{strictly quasi-pure resolution}
if $m_i>M_{i-1}$ for all $i$.
In this case we show that the bounds in (\ref{guess2}) are valid.
The idea of the proof is similar to the one  of \cite[Theorem 1.2]{HESR98}.

\begin{thm}
\label{quasipure}
Let $I \subset S$ be a graded ideal such that $R=S/I$
is Cohen--Macaulay which has a strictly quasi-pure resolution
and let $p=\projdim(R)$.
Then:
\begin{enumerate}
 \item
We have for $i=1,\dots,p$ that
$$
\beta_i^S(R)
\leq
\prod_{1\leq j<i}
\frac{M_j}{m_i-M_j}
\prod_{i<j\leq p}
\frac{M_j}{m_j-M_i}
\leq
\frac{1}{(i-1)!\cdot (p-i)!}
\prod_{j\neq i}  M_j.
$$
The first upper bound is reached for all $i$ if and only if $R$ has a pure resolution.
Every upper bound is reached for all $i$ if and only if $I$ has a linear resolution.
\item
We have for $i=1,\dots,p$ that
$$
\beta_i^S(R)
\geq
\prod_{1\leq j< i}
\frac{m_j}{M_i-m_j}
\cdot
\prod_{i<j\leq p}
\frac{m_j}{M_j-m_i}.
$$
Every lower bound is reached for all $i$ if and only if
$R$ has a pure resolution.
\end{enumerate}
\end{thm}
\begin{proof}
We consider the $(p\times p)$-square matrix
$$
A
=
\left(
  \begin{array}{cccc}
  \sum_{j} \beta_{1j}^S(R)                      &
\sum_{j} \beta_{2j}^S(R)                           &
\cdots & \sum_{j} \beta_{pj}^S(R) \\
  \sum_{j} j \cdot \beta_{1j}^S(R) &
  \sum_{j} j \cdot \beta_{2j}^S(R)       & \cdots &
  \sum_{j} j \cdot \beta_{pj}^S(R)   \\
  \vdots                            & \vdots                            & \vdots & \vdots \\
  \sum_{j} j^{p-1} \cdot \beta_{1j}^S(R) &
  \sum_{j} j^{p-1} \cdot \beta_{2j}^S(R)       & \cdots &
  \sum_{j} j^{p-1} \cdot \beta_{pj}^S(R)   \\
  \end{array}
\right)
.$$
We compute the determinant of $A$ as
$$
\det(A)
=
\sum_{j_1}
\dots
\sum_{j_p}
V(j_1,\dots,j_p)
\cdot
\prod_{1\leq l \leq p} \beta_{lj_l}^S(R)
$$
with the Vandermonde determinants
$$
V(j_1,\dots,j_p)
=
\det
\left(
  \begin{array}{ccccc}
    1           & 1             & \cdots    & 1 \\
    j_1     & j_2       & \cdots    & j_p   \\
    \vdots      & \vdots    & \vdots    & \vdots \\
    j_1^{p-1}   & j_2^{p-1}     & \cdots    & j_p^{p-1}   \\
  \end{array}
\right)
.
$$
Since $R$ has a strictly quasi-pure resolution
we have that $j_i > j_k$ for all integers $i,k$ such that
$i>k$, $\beta_{ij_i}^S(R)\neq 0$ and
$\beta_{kj_k}^S(R)\neq 0$.
Thus all the involved Vandermonde determinants are always positive.

We may compute $\det(A)$ also in a different way.
Fix $i \in \{1,\dots,p\}$.
By replacing the $i$-th column of $A$ by the alternating sum of all columns of $A$,
we obtain a matrix $A'$ such that $ \det (A) = \det(A')$.
It follows from Lemma \ref{quasihelpers},
that the $i$-th column of $A'$ is the transpose of the vector
$((-1)^{i+1},0,\dots,0)$.
Hence by expanding the determinant of $A'$ with respect to
the $i$-th column, we get
$$
\det(A)=\det(A')= \det(B),
$$
where $B$ is the $(p-1\times p-1)$-matrix
$$
\left(
  \begin{array}{cccccc}
    \sum_{j} j \beta_{1j}^S(R)      &\cdots&
    \sum_{j} j\beta_{i-1j}^S(R) &
    \sum_{j} j \beta_{i+1j}^S(R)   &\cdots&
    \sum_{j} j \beta_{pj}^S(R)\\
    \vdots & \vdots &  \vdots &  \vdots & \vdots    & \vdots \\
    \sum_{j} j^{p-1}  \beta_{1j}^S(R)   &\cdots&
    \sum_{j} j^{p-1} \beta_{i-1j}^S(R)  &
    \sum_{j} j^{p-1}\beta_{i+1j}^S(R)   &\cdots&
    \sum_{j} j^{p-1} \beta_{pj}^S(R)\\
  \end{array}
\right).
$$
Thus
\begin{eqnarray*}
\det(A) &=& \det(B)\\
&=&
\sum_{j_1}
\dots
\sum_{j_{i-1}}
\sum_{j_{i+1}}
\dots
\sum_{j_p}
U(j_1, \dots,j_{i-1},j_{i+1},\dots,j_p)
\cdot
\prod_{1\leq l \leq p, l\neq i} j_l \cdot \beta_{lj_l}^S(R)
\end{eqnarray*}
with the corresponding Vandermonde determinants
$$
U(j_1, \dots,j_{i-1},j_{i+1},\dots, j_p)
$$
$$
=
\det
\left(
  \begin{array}{cccccc}
    1              & \cdots & 1                    & 1                    & \cdots  & 1 \\
    j_1            & \cdots & j_{i-1}              & j_{i+1}              & \cdots  & j_p   \\
    \vdots         & \vdots & \vdots               & \vdots               & \vdots  & \vdots\\
    j_1^{p-2}      & \cdots & j_{i-1}^{p-2}        & j_{i+1}^{p-2}        & \cdots  & j_p^{p-2}   \\
  \end{array}
\right).
$$
Observe that
$$
V(j_1,\dots,j_p)
=
\prod_{i<l\leq p} (j_l - j_i)
\cdot
\prod_{1\leq l<i} (j_i - j_l)
\cdot
U(j_1, \dots,j_{i-1},j_{i+1},\dots, j_p).
$$
All in all we obtain from the discussion so far that
\begin{eqnarray}
\label{ende}
&&
\sum_{j_1}
\dots
\sum_{j_p}
\prod_{i<l\leq p} (j_l - j_i)
\cdot
\prod_{1\leq l<i} (j_i - j_l)
\cdot
U(j_1, \dots,j_{i-1},j_{i+1},\dots, j_p)
\cdot
\prod_{1\leq l \leq p} \beta_{lj_l}^S(R)
\end{eqnarray}
$$
=
\sum_{j_1}
\dots
\sum_{j_{i-1}}
\sum_{j_{i+1}}
\dots
\sum_{j_p}
U(j_1, \dots,j_{i-1},j_{i+1},\dots,j_p)
\cdot
\prod_{1\leq l \leq p, l\neq i} j_l
\cdot
\prod_{1\leq l \leq p, l\neq i}
\beta_{lj_l}^S(R).
$$
It follows from the fact that $R$ has a strict
quasi-pure resolution that for all integers $l$
\begin{eqnarray}
\label{twomore1}
\prod_{1\leq l<i} (m_i - M_l)
\prod_{i<l\leq p} (m_l - M_i)
\end{eqnarray}
$$
\leq
\prod_{1\leq l<i} (j_i - j_l)
\prod_{i<l\leq p} (j_l - j_i)\leq
\prod_{1\leq l<i} (M_i - m_l)
\prod_{i<l\leq p} (M_l - m_i).
$$
We always have for those $j_l$ with $\beta^S_{lj_l} (R) \neq 0$
that
\begin{eqnarray}
\label{twomore2}
\prod_{1\leq l \leq p, l\neq i} m_l
\leq
\prod_{1\leq l \leq p, l\neq i} j_l
\leq
\prod_{1\leq l \leq p, l\neq i} M_l.
\end{eqnarray}
Using (\ref{ende}),
the lower bound of (\ref{twomore1})
and the upper bound of (\ref{twomore2})
we obtain
$$
\beta_i^S(R)
\cdot
\prod_{i<l\leq p} (m_l - M_i)
\prod_{1\leq l<i} (m_i - M_l)
$$
$$
\cdot
\sum_{j_1}
\dots
\sum_{j_{i-1}}
\sum_{j_{i+1}}
\dots
\sum_{j_p}
U(j_1, \dots,j_{i-1},j_{i+1},\dots, j_p)
\cdot
\prod_{1\leq l \leq p, l\neq i}
\beta_{lj_l}^S(R)
$$
$$
\leq
\prod_{1\leq l \leq p, l\neq i} M_l
\sum_{j_1}
\dots
\sum_{j_{i-1}}
\sum_{j_{i+1}}
\dots
\sum_{j_p}
U(j_1, \dots,j_{i-1},j_{i+1},\dots, j_p)
\cdot
\prod_{1\leq l \leq p, l\neq i}
\beta_{lj_l}^S(R)
$$
and thus
$$
\beta_i^S(R)
\leq
\prod_{1\leq j < i} \frac{M_j}{m_i - M_j}
\cdot
\prod_{i<j\leq p} \frac{M_j}{m_j - M_i}.
$$
Analogously
using (\ref{ende}),
the upper bound of (\ref{twomore1})
and the lower bound of (\ref{twomore2})
we get
$$
\beta_i^S(R)
\geq
\prod_{1\leq j < i} \frac{m_j}{M_i - m_j}
\cdot
\prod_{i<j\leq p} \frac{m_j}{M_j - m_i}.
$$
(These lower bounds will also be proved in Section \ref{boijsoederberg} where we show that they hold in general for all Cohen--Macaulay rings.)
Checking the inequalities we see that
we have equality above or below for all  $1\leq i \leq p$
if and only if $R$ has a pure resolution.
We already observed that we have the inequalities $\prod_{1\leq j<i}
\frac{M_j}{m_i-M_j}
\prod_{i<j\leq p}
\frac{M_j}{m_j-M_i}
\leq
\frac{1}{(i-1)!\cdot (p-i)!}
\prod_{j\neq i}  M_j$. A straightforward discussion shows that
$\beta_i^S(R)=\frac{1}{(i-1)!\cdot (p-i)!}
\prod_{j\neq i}  M_j$ for all $i$ if and only if $I$ has a linear resolution.
This concludes the proof.
\end{proof}


\section{Componentwise linear ideals}
\label{componentwise}
Let $I \subset S=K[x_1,\dots,x_n]$ be a graded ideal.
Recall that $I$ has a \emph{$k$-linear resolution} if $\beta_{i,i+j}^S(I)=0$ for $j\neq k$.
For a non-negative integer $k$
we denote by $I_{\langle k \rangle} \subset S$ the ideal which is generated by all elements in $I_k$.
Herzog and Hibi \cite{HEHI99}
called $I$ \emph{componentwise linear} if $I_{\langle k \rangle}$ has a $k$-linear resolution for all $k\geq 0$.

It is well-known that a lot of important classes of ideals in combinatorial commutative algebra
are componentwise linear.
Recall that an ideal $I \subset S$ is called a \emph{monomial ideal} if it is
generated by monomials of $S$. Then we denote by $G(I)$
the unique minimal system of generators of $I$.
A monomial ideal $I \subset S$ is called \emph{strongly stable},
if  for all monomials $x^u=\prod_{k=1}^n x_k^{u_k} \in G(I)$
and $i$ with $x_i|x^u$ we have for all $1\leq j\leq i$ that $(x^u/x_i)x_j\in I$.
It is well-known that strongly stable ideals are componentwise linear.
But also stable ideals, squarefree (strongly) stable ideals and more generally
$\ab$-stable ideal  are componentwise linear. (See \cite[Theorem 3.11]{RO05} for definitions and a proof.)
In particular, this implies that all generic initial ideals are componentwise linear
provided $\chara(K)=0$.
(E.g.\ see \cite{ARHEHI00}  or \cite[Lemma 3.3]{HEZH}.)

In the proof of the next theorem we will need the
Eliahou--Kervaire formula \cite{ElKe} for the graded Betti-numbers of a strongly stable ideal $I$:
we have for all $i\geq 1$ and $j\geq 0$ that
$
\beta_{i,i+j}^S(S/I)
=
\sum_{x^u \in G(I),\ \deg x^u=j+1} \binom{m(u)-1}{i-1}
$
where we set $m(u)=\max\{i: 1\leq i \leq n,\ u_i>0\}$ for a monomial $x^u$ with $u \in \N^n$.
Here we make the convention that $\binom{a}{b}=0$ for $a,b \in \Z$ unless $0\leq b \leq a$.
(Note that these formulas above are already true for stable ideals.)
Observe the following facts.
If $\beta_{i,i+j}^S(S/I) \neq 0$ for some $i$, then
$\beta_{k,k+j}^S(S/I)\neq 0$ for $1\leq k \leq i$.
Moreover, only those $x^u \in G(I)$ with $m(u)\geq i$
are relevant for the total $i$-th Betti number $\beta_{i}^S(S/I) \neq 0$
and then $\deg x^u \geq m_i-i+1$.

\begin{thm}
\label{componentwise1}
Let $I \subset S$ be a componentwise linear ideal such that $R=S/I$ is Cohen--Macaulay
and let $p=\projdim(R)$.
Then:
\begin{enumerate}
 \item
We have for $i=1,\dots,p$ that
$$
\beta_i^S(R)
\leq
\prod_{1\leq j<i}
\frac{M_j}{M_i-M_j}
\cdot
\prod_{i<j\leq p}
\frac{M_j}{M_j-M_i}
\leq
\frac{1}{(i-1)!\cdot (p-i)!}
\prod_{j\neq i}  M_j.
$$
Every upper bound is reached for all $i$ if and only if $I$ has a linear resolution.
\item
We have for $i=1,\dots,p$ that
$$
\beta_i^S(R)
\geq
\prod_{1\leq j< i}
\frac{m_j}{m_i-m_j}
\cdot
\prod_{i<j\leq p}
\frac{m_j}{m_j-m_i}
\geq
\prod_{1\leq j< i}
\frac{m_j}{M_i-m_j}
\cdot
\prod_{i<j\leq p}
\frac{m_j}{M_j-m_i}.
$$
Every lower bound is reached for all $i$ if and only if $I$ has a linear resolution.
\end{enumerate}
\end{thm}
\begin{proof}
Without loss of generality we may assume that the field $K$ is infinite.
We denote by $\gin(I)$ the generic initial ideal of $I$ with respect to
the reverse lexicographical order.
The proof of the main result in \cite{ARHEHI00} and \cite[Lemma 3.3]{HEZH}
shows that $\gin(I)$ has the same graded Betti numbers as $I$ and is
a stable ideal in all characteristics.
If we replace $I$ by $\gin(I)$,
then the Betti numbers of $I$ do not depend on the characteristic
of $K$ and we may assume that $\chara(K)=0$.
Replacing another time $I$ by $\gin(I)$ does not change the Betti numbers and
thus we may now assume that $I$ is a strongly stable ideal.

Since $R$ is Cohen--Macaulay and it is known that
$x_n,\dots,x_{n-\depth(R)+1}$ is a regular sequence for $R$,
we may assume that $\dim(R)=0$ and
thus a pure power of each variable belongs to $I$.
Let $a>0$ be the smallest natural number such that $x_n^a \in I$.
Then $\deg(x^u) \leq a$ for all $x^u \in G(I)$
and $x_n^a \in G(I)$, because $I$ is strongly stable.
Note that for (i) and (ii) we have to show only the corresponding
first inequalities, since the other are trivially true as noted in the other sections
of this paper.

(i):
It follows from the Eliahou--Kervaire formula for the graded Betti numbers of $R$
that
$$
M_i= a + i-1  \text{ for } i=1,\dots,n.
$$
We have that $(x_1,\dots,x_n)^a \subseteq I$ and thus it follows from
\cite[Theorem 3.2]{COHEHI04}
that
$$
\beta^S_i(S/I)
\leq
\beta^S_i(S/(x_1,\dots,x_n)^a)
=
\prod_{1\leq j<i}
\frac{M_j}{M_i-M_j}
\cdot
\prod_{i<j\leq n}
\frac{M_j}{M_j-M_i}
$$
where the last equation follows from the fact that
$(x_1,\dots,x_n)^a$ has an $a$-linear resolution,
the maximal shifts coincide with the ones of $I$
and that in this case the equation follows from \cite[Theorem 1]{HEKU84}).

If $I$ has a linear resolution, then
$I=(x_1,\dots,x_n)^a$ and
the upper bounds for $\beta_i^S(S/I)$ are reached.
Assume that we have equations everywhere.
Then it follows that
$\beta_i^S(S/I)=\beta_i^S(S/(x_1,\dots,x_n)^a)$ for $i=1,\dots,n$.
In the proof of \cite[Theorem 3.2]{COHEHI04}
it is shown, that this implies
$$
|x^u \in G(I_{\langle j \rangle}) : m(u)=k |
=
|x^u \in G((x_1,\dots,x_n)^a_{\langle j \rangle}) : m(u)=k |
\text{ for } j \in \Z,\ 1\leq k \leq n-1.
$$
This implies that $I_{\langle j \rangle}=0$ for $j<a$
and thus $I=(x_1,\dots,x_n)^a$.
Hence $I$ has an $a$-linear resolution.

(ii):
Fix $1\leq i \leq n$ and write $m_i=i+b-1$ for some natural number $b$.
Let $J=I_{\geq b}$ be the ideal which is generated by all elements of $I$
of degree greater or equal to $b$.
It follows from the Eliahou--Kervaire formula
and the observations given above
that
\begin{eqnarray*}
\beta_j^S(S/J)&=&\beta_j^S(S/I) \text{ for } j \geq i,\\
m_j(S/J) &=& m_j(S/I) \text{ for } j \geq i,\\
m_j(S/J)&\geq & m_j(S/I) \text{ for } 1\leq j < i,\\
m_j(S/J) &=& m_i(S/J)-(i-j) \text{ for } 1\leq j < i.
\end{eqnarray*}
Note that $S/J$ is still zero dimensional.
Assume that we could prove the lower bound for $S/J$,
then it would follow that
\begin{eqnarray*}
\beta_i^S(S/I)&=& \beta_i^S(S/J)
\\
 &\geq &
\prod_{1\leq j< i}
\frac{m_j(S/J)}{m_i(S/J)-m_j(S/J)}
\cdot
\prod_{i<j\leq n}
\frac{m_j(S/J)}{m_j(S/J)-m_i(S/J)}
\\
&= &
\prod_{1\leq j< i}
\frac{m_j(S/J)}{m_i(S/I)-m_j(S/J)}
\cdot
\prod_{i<j\leq n}
\frac{m_j(S/I)}{m_j(S/I)-m_i(S/I)}
\\
&\geq &
\prod_{1\leq j< i}
\frac{m_j(S/I)}{m_i(S/I)-m_j(S/I)}
\cdot
\prod_{i<j\leq n}
\frac{m_j(S/I)}{m_j(S/I)-m_i(S/I)}.
\end{eqnarray*}
The last inequality follows because
for $1\leq j< i$ we have
\begin{eqnarray*}
&& \frac{m_j(S/J)}{m_i(S/I)-m_j(S/J)} \geq \frac{m_j(S/I)}{m_i(S/I)-m_j(S/I)}\\
&\Leftrightarrow &
m_j(S/J)m_i(S/I) - m_j(S/J) m_j(S/I) \geq m_j(S/I)m_i(S/I) -m_j(S/I)m_j(S/J)\\
&\Leftrightarrow &
m_j(S/J)m_i(S/I)  \geq m_j(S/I)m_i(S/I) \\
&\Leftrightarrow &
m_j(S/J) \geq m_j(S/I).
\end{eqnarray*}
Here the last inequality follows from the definition of $J$ as noted above.
It remains to show the lower bound for $\beta_i^S(S/J)$.
Let $L=(x_1,\dots,x_n)^b$.
We observe that $J\subseteq L$ and we have
\begin{eqnarray*}
m_j(S/L) &=& m_i(S/L)-(i-j) \text{ for } 1\leq j < i,\\
m_j(S/L) &=& m_i(S/L)+(j-i) \text{ for } i< j \leq n,\\
m_j(S/J) &=& m_j(S/L) \text{ for } j \leq i,\\
m_j(S/J) &\geq & m_j(S/L) \text{ for } i< j \leq n.\\
\end{eqnarray*}
Moreover, it follows  from \cite[Theorem 3.2]{COHEHI04}
that $\beta_i^S(S/J) \geq  \beta_i^S(S/L)$.
We compute
\begin{eqnarray*}
\beta_i^S(S/J)&\geq & \beta_i^S(S/(x_1,\dots,x_n)^b)
\\
&= &
\prod_{1\leq j< i}
\frac{m_j(S/L)}{m_i(S/L)-m_j(S/L)}
\cdot
\prod_{i<j\leq n}
\frac{m_j(S/L)}{m_j(S/L)-m_i(S/L)}
\\
&= &
\prod_{1\leq j< i}
\frac{m_j(S/J)}{m_i(S/J)-m_j(S/J)}
\cdot
\prod_{i<j\leq n}
\frac{m_j(S/L)}{m_j(S/L)-m_i(S/J)}
\\
&\geq &
\prod_{1\leq j< i}
\frac{m_j(S/J)}{m_i(S/J)-m_j(S/J)}
\cdot
\prod_{i<j\leq n}
\frac{m_j(S/J)}{m_j(S/J)-m_i(S/J)}
\end{eqnarray*}
The last inequality follows because for $i< j\leq n$ we have
\begin{eqnarray*}
&& \frac{m_j(S/L)}{m_j(S/L)-m_i(S/J)} \geq \frac{m_j(S/J)}{m_j(S/J)-m_i(S/J)}\\
&\Leftrightarrow &
m_j(S/L)m_j(S/J)-m_j(S/L)m_i(S/J) \geq m_j(S/J)m_j(S/L)- m_j(S/J)m_i(S/J)\\
&\Leftrightarrow &
-m_j(S/L)m_i(S/J) \geq - m_j(S/J)m_i(S/J)\\
&\Leftrightarrow &
m_j(S/J) \geq m_j(S/L).
\end{eqnarray*}
The last inequality is valid as noted above.
Thus we get the desired lower bound for $J$ and hence also for $I$.

Assume that for all $i$ the lower bound for $\beta_i^S(S/I)$ is reached.
For $i=1$ the corresponding constructed $J$ is just $I$.
It follows then also that $\beta_1^S(S/I)=\beta_1^S(S/L)$
and applying again \cite[Theorem 3.2]{COHEHI04} we see that
$\beta_i^S(S/I)=\beta_i^S(S/L)$ for $1\leq i \leq n$.
Now we deduce as in the proof of (i) that
indeed $I$ has a linear resolution.
This concludes the proof.
\end{proof}

The Cohen--Macaulay assumption is essential for the lower bounds in (\ref{guess2}).
In fact, we can construct a strongly stable ideal as a counterexample.
The ideal is taken from \cite{HESR98}.

\begin{ex}
\label{counterguess2}
Let $S=K[x_1,\dots,x_4]$ be a polynomial ring in 4 variables
and we consider the strongly stable ideal $I=(x_1^2,x_1x_2, x_2^3, x_2^2x_3, x_2^2x_4)$.
Then $S/I$ is not Cohen--Macaulay because $\dim(S/I)=2$ and $\depth(S/I)=0$.
It follows from the Eliahou--Kervaire formula that
$$
\beta_1^S(S/I)=5,\
\beta_2^S(S/I)=7,\
\beta_3^S(S/I)=4,\
\beta_4^S(S/I)=1
$$
and
$$
m_1=2,\
m_2=3,\
m_3=5,\
M_4=m_4=6.
$$
But now
$$
\beta_4^S(S/I)=1 <
\frac{m_3}{M_4-m_3} \cdot \frac{m_2}{M_4-m_2} \cdot \frac{m_1}{M_4-m_1}
=
\frac{5\cdot 3\cdot 2}{1 \cdot 3 \cdot 4}
=
\frac{30}{12}.
$$
\end{ex}

On the other hand for strongly stable ideals we still can give
an upper bound for the $i$-th total Betti number without the Cohen--Macaulay assumption.

\begin{thm}
\label{componentwise2}
Let $I \subset S$ be a componentwise linear ideal and $p=\codim(S/I)$.
We have for $i=1,\dots,p$ that
$$
\beta_i^S(S/I)
\leq
\binom{i+M_1-2 }{i-1} \cdot \binom{p+M_1-1}{p-i}.
$$
The upper bound is reached for all $i$ if and only if
$S/I$ is Cohen--Macaulay and $I$ has a linear resolution.
\end{thm}
\begin{proof}
As shown in the proof of Theorem \ref{componentwise1}  we may assume that $\chara(K)=0$ and
that $I$ is a strongly stable ideal.
It is known that $x_n,\dots,x_{n-\depth(S/I)+1}$ is a regular sequence for $S/I$ and
thus we may assume that $\depth(S/I)=0$, i.e. $\projdim (S/I)=n$.

Let $J=I_{\geq M_1(S/I)}$ be the ideal which is generated by all elements of $I$
of degree greater or equal to $M_1(S/I)$.
It follows  from \cite[Theorem 3.2]{COHEHI04} that $\beta_i^S(S/I) \leq  \beta_i^S(S/J)$.
Note that $ M_i\leq M_1 + i-1$ as one deduces from the
Eliahou--Kervaire formula.
By construction of $J$ we have $M_1(S/J)=M_1$ and
$J$ has an $ M_1$-linear resolution.
Assuming that we can show the upper bound for $S/J$,
we get for $1\leq i \leq n$ that
$$
\beta_i^S(S/I)\leq  \beta_i^S(S/J)\leq
\binom{i+M_1-2 }{i-1} \cdot \binom{n+M_1-1}{n-i}.
$$
It remains to show the upper bound for $S/J$.
Note that
\begin{eqnarray*}
\beta_i^S(S/J)
&=&
\sum_{x^u \in G(J)} \binom{m(u)-1}{i-1}\\
&=&
\sum_{j=i}^n \binom{j-1}{i-1} |\{x^u \in G(J) : m(u)=j  \}|\\
&\leq&
\sum_{j=i}^n
\binom{j-1}{i-1} \binom{j+M_1-1-1}{M_1-1}.
\end{eqnarray*}
We prove by induction on $n-i$
that
$$
\sum_{j=i}^n
\binom{j-1}{i-1} \binom{j+M_1-2}{M_1-1}
=
\binom{i+M_1-2 }{i-1} \cdot \binom{n+M_1-1}{n-i}.
$$
The assertion is trivial for $i=n$. Let $i<n$.
Using the induction hypothesis we compute
\begin{eqnarray*}
&&
\sum_{j=i}^n \binom{j-1}{i-1} \binom{j+M_1-2}{j-1}
\\
&=&
\sum_{j=i}^{n-1} \binom{j-1}{i-1} \binom{j+M_1-2}{j-1} +  \binom{n-1}{i-1} \binom{n+M_1-2}{n-1}
\\
&=&
\binom{i+M_1-2 }{i-1} \cdot \binom{n-1+M_1-1}{n-1-i}    +  \binom{n-1}{i-1} \binom{n+M_1-2}{n-1}
\\
&=&
\binom{i+M_1-2 }{i-1} \cdot \binom{n+M_1-1}{n-i}
\\
&&
-\binom{i+M_1-2 }{i-1}\binom{n+M_1-2}{n-i}
+  \binom{n-1}{i-1} \binom{n+M_1-2}{n-1}
\\
&=&
\binom{i+M_1-2 }{i-1} \cdot \binom{n+M_1-1}{n-i}
\end{eqnarray*}
because
$\binom{i+M_1-2 }{i-1}\binom{n+M_1-2}{n-i}=\binom{n-1}{i-1} \binom{n+M_1-2}{n-1}$
as one verifies by a direct computation.
If $S/I$ is Cohen--Macaulay and $I$ has a ($M_1$-)linear resolution,
then we know by \cite{HEKU84} that $\beta_i^S(S/I)$ reaches the upper bound
for all $i$. Assume now that
$
\beta_i^S(S/I)
=
\binom{i+M_1-2 }{i-1} \cdot \binom{p+M_1-1}{p-i}$
for $i=1,\dots,p$. As seen above we may assume that $p=n$.
Then the corresponding bounds for $S/J$ are also achieved.
It follows from the inequalities above that
 every monomial of degree $M_1$ is a minimal generator of $J$.
This means that $J=(x_1,\dots,x_n)^{M_1}$.
Thus $S/J$ is zero dimensional and hence Cohen--Macaulay.
But since $J=I_{\geq M_1}$ then also $S/I$ is zero dimensional and therefore Cohen--Macaulay.
Now we can apply Theorem \ref{componentwise1} (i) to conclude that $I$ has a linear resolution.
\end{proof}

\begin{rem}
The results of this section can also be used to prove bounds for the Betti numbers if
$I$ is not componentwise linear, at least if $\chara(K)=0$. Let $I \subset S$ be an arbitrary graded ideal and $p=\projdim(S/I)$. Recall that  $\reg(S/I)=\max_{1\leq i \leq \projdim(S/I)} \{ M_i -i\}$ is called the \emph{Castelnuovo--Mumford regularity} of $S/I$. It is well-known that
$\reg(S/I)=\reg(S/\gin(I))$ where $\gin(I)$ is the generic initial ideal of $I$ with respect to
the reverse lexicographical order (see \cite{BAST}). Moreover, $\beta_i^S(S/I)\leq \beta_i^S(S/\gin(I))$ for all $i$. Since $\gin(I)$ is componentwise linear it follows from these observations and Theorem \ref{componentwise2} that
\begin{eqnarray*}
\beta_i^S(S/I)
&\leq&
\beta_i^S(S/\gin(I))\\
&\leq&
\binom{i+\reg(S/\gin(I))-1 }{i-1} \cdot \binom{p+\reg(S/\gin(I))}{p-i}\\
&=&
\binom{i+\reg(S/I)-1 }{i-1} \cdot \binom{p+\reg(S/I)}{p-i}
\end{eqnarray*}
where $p=\codim(S/I)=\codim(S/\gin(I))$. With similar arguments one can use Theorem \ref{componentwise1} to prove upper bounds in the Cohen--Macaulay case using the regularity. Since we get better results for this case in the next section, we omit the details.
\end{rem}

\section{Cohen--Macaulay rings}
\label{boijsoederberg}

We saw that the lower and upper bounds in (\ref{guess1}) do not hold in general. Also
the upper bounds in (\ref{guess2}) are not candidates for upper bounds since the numbers may be negative. Using the Boij--S\"oderberg theory which was conjectured and developed partly by Boij--S\"oderberg  \cite{BOSO} and then completely  proved by Eisenbud-Schreyer \cite{EISCH} (see also \cite{BOSO} and \cite{EIFLWE}) we show that
the lower bounds in (\ref{guess2}) and upper bounds in (\ref{guess3}) hold under the Cohen--Macaulay assumption.

We recall parts of the Boij--S\"oderberg theory which is needed in the following. Fix a positive integer $p$.
For any strictly increasing sequence of non-negative integers $d=(d_0,d_1,\dots,d_p)$ with $d_0=0$
we define a diagram $\pi(d)$ by
$$
\pi(d)_{i,j}
=
\begin{cases}
\prod_{1\leq k <i} \frac{d_k}{d_i-d_k}\prod_{i< k \leq p} \frac{d_k}{d_k-d_i} & \text{if } j=d_i,\\
0 & \text{else}
\end{cases}
$$
and call $\pi(d)$ a \emph{pure diagram}. The sequence $d=(d_0,d_1,\dots,d_p)$ is called the \emph{degree sequence} of the diagram.  Note that there is a choice which diagrams are called the pure ones up to multiplication with respect to a positive real number. We choose them in such a way that $\pi(d)_{0,0}=1$. A pure diagram is called \emph{linear} if $d_k=d_1+(k-1)$ for $1\leq k \leq p$.
There exists a partial order on pure diagrams by defining
$\pi(d)\leq \pi(d')$ for two
increasing sequences of non-negative integers $d=(0,d_1,\dots,d_p)$
and $d'=(0,d'_1,\dots,d'_p)$ if and only if $d\leq d'$ coefficientwise.
For two fixed increasing sequences of positive integers $\underline{d}$ and
$\overline{d}$ denote by $\Pi_{\underline{d},\overline{d}}$
the set of  pure diagrams $\pi(d)$
such that $\pi(\underline{d})\leq \pi(d)\leq \pi(\overline{d})$.
Since $\pi(d)_{i,j}\neq 0$ only for finitely many $i,j$
we can
consider the convex hull of  $\Pi_{\underline{d},\overline{d}}$,
that is the set of convex combinations
$
D=\sum_{\pi(d) \in \Pi_{\underline{d},\overline{d}}} \lambda_d \pi(d)
$
with real non-negative coefficients $\lambda_d$
and
$\sum_{\pi(d) \in \Pi_{\underline{d},\overline{d}}} \lambda_d =1$.

One of the main results of the Boij--S\"oderberg theory implies (see \cite[Conjecture 2.4]{BOSO} and the full proof in \cite{EISCH}) that for a Cohen--Macaulay algebra $R$ of projective dimension $p$
the Betti-diagram $\beta^S(R)=(\beta^S_{i,j}(R))$
is a convex combination of the convex hull of  $\Pi_{m,M}$
where $m=(m_1,\dots,m_p)$ and $M=(M_1,\dots,M_p)$
and the $m_i,M_i$ are the usual maximal and minimal shifts
in the minimal graded free resolution of $R$.

Note that the Boij--S\"oderberg theory treads
more generally modules instead of rings. Then one of the results is that the Betti diagram of a Cohen--Macaulay module may be written (uniquely) as a positive rational linear combination of pure diagrams whose degree sequences form a totally ordered sequence. Since $\beta^S_{0,0}(R)=1$ and $\beta^S_{0,j}(R)=0$ for $j\neq 0$, the Betti diagram of $R$ is already a convex combination of pure diagrams as considered above and we restrict ourself to this situation.

Now we consider the  convex hull of  $\Pi_{\underline{d},\overline{d}}$,
and a convex combination $D$ as described above.
We define formally for $0\leq i \leq p$ and $j \in \Z$ the numbers
$$
\beta_{i,j}(D)=
\sum_{\pi(d) \in \Pi_{\underline{d},\overline{d}}}
\lambda_d \pi(d)_{i,j},\quad
\beta_{i}(D)=\sum_{j\in \Z} \beta_{i,j}(D).
$$
We also set for $1\leq i \leq p$
$$
M_i(D)=\max\{j \in \Z \colon \beta_{i,j}(D)\neq 0\}
\text{ and }
m_i(D)=\min\{j \in \Z \colon \beta_{i,j}(D)\neq 0\}.
$$
Observe that
$$
M_i(D)=\max\{ d_i : \lambda_d\neq 0 \}
\text{ and }
m_i(D)=\min\{ d_i : \lambda_d \neq 0 \}.
$$
Note also that it follows from the definition of the diagrams $\pi(d)$
that $M_i(D)<M_{i+1}(D)$ and $m_i(D)<m_{i+1}(D)$ hold for $1\leq i<p$.
At first we prove the following purely numerical result.
\begin{thm}
\label{thm:cmideals}
Let
$\underline{d}=(\underline{d}_0,\dots,\underline{d}_p)$
and
$\overline{d}=(\overline{d}_0,\dots,\overline{d}_p)$
be two strictly increasing sequences of non-negative integers
with $\underline{d}_0=\overline{d}_0=0$
such that
$\underline{d} \leq \overline{d}$.
Assume that
$D=\sum_{\pi(d) \in \Pi_{\underline{d},\overline{d}}} \lambda_d \pi(d)$
is a convex combination of elements of $\Pi_{\underline{d},\overline{d}}$.
Then:
\begin{enumerate}
 \item
We have for $i=1,\dots,p$ that
$$
\beta_i(D)
\leq
\frac{1}{(i-1)!\cdot (p-i)!}
\prod_{j\neq i}  M_j(D).
$$
The upper bound is reached for all $i$
if and only if  $D$ is a linear diagram.
\item
We have for $i=1,\dots,p$ that
$$
\beta_i(D)
\geq
\prod_{1\leq j< i}
\frac{m_j(D)}{M_i(D)-m_j(D)}
\cdot
\prod_{i<j\leq p}
\frac{m_j(D)}{M_j(D)-m_i(D)}.
$$
Every lower bound is reached for all $i$
if and only if  $D$ is a pure diagram.
\end{enumerate}
\end{thm}
\begin{proof}
(i)
We compute
\begin{eqnarray*}
\beta_{i}(D)
&=&
\sum_{\pi(d) \in \Pi_{\underline{d},\overline{d}}}
\lambda_d \pi(d)_{i}\\
&=&
\sum_{\pi(d) \in \Pi_{\underline{d},\overline{d}}}
\lambda_d \prod_{1\leq j <i} \frac{d_j}{d_i-d_j}\cdot\prod_{i< j \leq p} \frac{d_j}{d_j-d_i}\\
&\leq&
\sum_{\pi(d) \in \Pi_{\underline{d},\overline{d}}, \lambda_d\neq 0  }
\lambda_d \prod_{1\leq j <i} \frac{M_j(D)}{i-j}\cdot\prod_{i< j \leq p} \frac{M_j(D)}{j-i}\\
&=&
\frac{1}{(i-1)!\cdot (p-i)!}
\prod_{j\neq i}  M_j(D).
\end{eqnarray*}
Note that if $D$ is not a pure diagram, then the inequality is strict.
But even for a pure diagram which is not linear the inequality is strict.
Hence we have equality if and only if $D$ is a linear diagram.

(ii)
Observe that $M_i(D)-m_j(D)\geq m_i(D)-m_j(D)>0$ for $j<i$
and similar $M_j(D)-m_i(D)>0$ for $i<j$.
Then we get
\begin{eqnarray*}
\beta_{i}(D)
&=&
\sum_{\pi(d) \in \Pi_{\underline{d},\overline{d}}}
\lambda_d \pi(d)_{i}\\
&=&
\sum_{\pi(d) \in \Pi_{\underline{d},\overline{d}}}
\lambda_d \prod_{1\leq j <i} \frac{d_j}{d_i-d_j}\cdot \prod_{i< j \leq p} \frac{d_j}{d_j-d_i}\\
&\geq&
\sum_{\pi(d) \in \Pi_{\underline{d},\overline{d}}, \lambda_d\neq 0}
\lambda_d \prod_{1\leq j <i} \frac{m_j(D)}{M_i(D)-m_j(D)} \cdot\prod_{i< j \leq p} \frac{m_j(D)}{M_j(D)-m_i(D)}.
\\
&=&
\prod_{1\leq j <i}\frac{m_j(D)}{M_i(D)-m_j(D)} \cdot\prod_{i< j \leq p} \frac{m_j(D)}{M_j(D)-m_i(D)}.
\end{eqnarray*}
Note that if $D$ is not a pure diagram, then the inequalities are strict in general.
Hence we have equalities for all $i$ if and only if $D$ is a pure diagram.
\end{proof}

As always $K$ is a field and $S=K[x_1,\dots,x_n]$ a standard graded polynomial ring.
As a corollary of Theorem \ref{thm:cmideals} and the Boij--S\"oderberg theory we get:
\begin{cor}
Let $I \subset S$ be a graded ideal such that $R=S/I$ is Cohen--Macaulay and let $p=\projdim(R)$.
Then:
\begin{enumerate}
 \item
We have for $i=1,\dots,p$ that
$$
\beta_i^S(R)
\leq
\frac{1}{(i-1)!\cdot (p-i)!}
\prod_{j\neq i}  M_j.
$$
The upper bound is reached for all $i$
if and only if  $I$ has a linear resolution.
\item
We have for $i=1,\dots,p$ that
$$
\beta_i^S(R)
\geq
\prod_{1\leq j< i}
\frac{m_j}{M_i-m_j}
\cdot
\prod_{i<j\leq p}
\frac{m_j}{M_j-m_i}.
$$
Every lower bound is reached for all $i$
if and only if  $R$ has a pure resolution.
\end{enumerate}
\end{cor}

\begin{rem}
Note that it is also known that the Betti diagram of a graded ring $S/I$ which is not necessarily
Cohen--Macaulay may be written (uniquely) as a positive rational linear combination of pure diagrams (see \cite{BOSO2}). But here the appearing degree sequences maybe of different lengths and this causes problems. Indeed the Cohen--Macaulay assumption is essential for the lower bound (\ref{guess2}) as we saw in Example (\ref{counterguess2}).
Similar upper bounds as the ones in (\ref{guess3}) can be proved in the case where $S/I$ is not Cohen--Macaulay. Since the formulas in this case are not as nice and compact as the ones in the Cohen--Macaulay case we do not present them here and leave the details to the interested reader.
\end{rem}

\end{document}